\documentstyle[12pt,amssymbols]{article}
\newcommand{\bq}{\begin{equation}}
\newcommand{\eq}{\end{equation}}
\newcommand{\bqa}{\begin{eqnarray}}
\newcommand{\eqa}{\end{eqnarray}}
\newcommand{\ra}{\rightarrow}

\def\g{\gamma}

\def\ep{\epsilon}

\def\ov{\over}

\def\ed{\end{document}}

\def\ra{\rightarrow}
\def\al{\alpha}
\def\2pi{1\over 2\pi i}
\def\q{q-q^{-1}}
\def\~{\tilde}
\def\newline{\hfil\break}

\def\ra{\rightarrow}

\def\sq2{\sqrt{2}}
\def\sqk2{\sqrt{2(k+2}}
\def\sqk{\sqrt{k}}

\def\be{\begin{equation}}
\def\ee{\end{equation}}
\def\br{\begin{array}}
\def\er{\end{array}}
\def\bea{\begin{eqnarray}}
\def\eea{\end{eqnarray}}
\def\ba{\begin{equation}\begin{array}}
\def\ea{\end{array}\end{equation}}
\def\bac{\begin{equation}\begin{array}{rll}}

\def\al{\alpha}

\newcommand{\LS}{\widehat{sl(2)}}
\newcommand{\uq}{U_q (\widehat{sl(2)})}

\newcommand{\Zep}{Z(\ep_{1},\dots,
\ep_{s}|z_{1},\dots,z_{s})}

\def\Z{{\Bbb Z}}
\def\C{{\Bbb C}}
\def\N{{\Bbb N}}

\def\ep{\epsilon}
\def\epp{\epsilon^{\prime}}
\def\Sh{S(\hat{h^{-}})}

\begin{document}
\rightline{CRM-2274}
\rightline{April 18, 1995}
\vbox{\vspace{-10mm}}
\vspace{1.0truecm}
\begin{center}
{\LARGE \bf A Quantum Analogue of the ${\cal Z}$ Algebra }\\[8mm]
{\large A. Hamid Bougourzi$^1$ and Luc Vinet$^2$ }\\
[6mm]{\it Centre de Recherches Math\'ematiques\\
Universit\'e de Montr\'eal\\
C.P. 6128-A, Montr\'eal, P.Q., H3C 3J7, Canada.}\\[20mm]
\end{center}
\vspace{1.0truecm}
\begin{abstract}
We define a natural quantum analogue for the ${\cal Z}$ algebra,
and which we refer to as the ${\cal Z}_q$ algebra, by modding out
the Heisenberg algebra from the
quantum affine algebra $U_q(\hat{sl(2)})$ with level
$k$.  We discuss the representation
theory of this ${\cal Z}_q$ algebra. In particular,
we exhibit its reduction to
a group
algebra, and to a tensor product of a group algebra with a quantum
Clifford algebra when $k=1$, and $k=2$, and thus,
we recover the explicit constructions of $\uq$-standard modules
as achieved by Frenkel-Jing and Bernard, respectively. Moreover,
for arbitrary nonzero level $k$, we show that the explicit basis for the
simplest ${\cal Z}$-generalized Verma module as constructed by
Lepowsky and primc is also a basis for its corresponding
${\cal Z}_q$-module, i.e., it is invariant under the $q$-deformation for
generic $q$.  We expect this ${\cal Z}_q$ algebra
(associated with $\uq$ at level $k$),
to play the role of
a dynamical symmetry in the off-critical $ Z_k$ statistical models.
\end{abstract}

\footnotetext[1]{
Email: {\tt bougourz@ere.umontreal.ca}}
\footnotetext[2]{
Email: {\tt vinet@ere.umontreal.ca}}
\section{Introduction}

One of the major recent developments in the field of integrable
models has been the realization by the Kyoto school
\cite{Daval92,JiMi94} of the important role played by non-Abelian and
dynamical symmetries in the resolution of integrable systems.
Prior to this and
besides conformal field theory, the main approach in the analysis of
integrable models  has been based on Abelian symmetries together with the
Bethe ansatz. However, this approach, despite its success of
being more systematic  in handling the spectra of
most integrable systems,  has its limitations as far as concrete
computations of physical quantities are concerned such as
correlation functions and form factors. The reason is that the
latter quantities are based on scalar products of the
eigenvectors of the Hamiltonians or transfer matrices of the
systems; but the eigenspaces in the physically interesting
thermodynamic limit are  infinite-dimensional, and hence it is
not easy to define their structures, and much less the
scalar products on them. However,
since some non-Abelian infinite-dimensional algebras have well
defined scalar products on their
infinite-dimensional modules, then if we succeed in
establishing that the Hamiltonian or the transfer matrix
of an integrable system commutes
with one of these algebras, we automatically know, not
only its eigenspaces  which are the modules of this algebra, but
also the scalar product on them. We might even describe
the local operators, and the
creation and annihilation operators of the eigenvectors in terms
of some operators related to this algebra, such as the intertwiners
of its modules, and hence we might be able to compute exactly the
correlation functions and the form factors.
In fact, this is precisely the program that has been
behind the
enormous success in the resolution of conformal
field theories, and more recently in the resolution of the
XXZ quantum spin chain model (which is equivalent to the
6-vertex classical model) in the antiferromagnetic regime by the
Kyoto school. We should however, mention Ref. \cite{Koral93}
where another approach to the calculation of correlation functions is
developed.

It is then an interesting program to build as many
infinite-dimensional algebras as possible, hoping that one of them
might turn out to be a non-Abelian or a dynamical symmetry
of an integrable model, and vice versa. The main point of this
paper is precisely to define a new infinite-dimensional
algebra, which, as explained in the next paragraph,
 should be the dynamical symmetry of the
 off-critical $Z_k$ models.

It is known that the $\LS$ affine Lie algebra is  a dynamical symmetry
in conformal
field theory (the continuum critical limit of the XXX model)
\cite{Aff88} and a
non-Abelian symmetry of the antiferromagnetic (off-critical) XXX model
\cite{JiMi94}. Moreover,  its deformation, the $\uq$ algebra, is also a
non-Abelian symmetry
of the
antiferromagnetic (off critical) XXZ model. It is also known
that  the critical $Z_k$ models, such as the Ising model ($k=2$) and the
Potts model ($k=3$),  have as a dynamical symmetry
the parafermionic algebra \cite{ZaFa84},
which is ironically related to the
${\cal Z}$ algebra \cite{LeWi84}. The latter algebra, in turn,
is obtained from
the quotient of the $\LS$ algebra with level $k$ by its Heisenberg subalgebra.

{}From all the above known results,
it is therefore natural to consider a similar construction
for a quantum analogue of the ${\cal Z}$ algebra, denoted by ${\cal Z}_q$,
 from the $\uq$ algebra, and
to expect it to play the same role of a dynamical symmetry for the
off-critical $ Z_k$ models.
In fact, such a program is already and implicitly
implemented in the simplest case of the off-critical Ising model ($k=2$),
where the ${\cal Z}_q$ algebra (to be precise the corresponding
quantum parafermionic algebra) reduces simply to
 a quantum Clifford algebra \cite{Foda94}.

This paper  is organized as follows: in section 2, we recall
basic definitions  about the $\uq$ quantum affine algebra, using
the formal variable approach. In section 3, we gradually introduce
the ${\cal Z}_q$
algebra by modding out the Heisenberg
subalgebra from $\uq$ with level $k$. We derive two defining relations,
called ``the quantum generalized commutation relations", for this
algebra, as well as the relations between its elements and
those of
$\uq$. In section 4, we discuss the reduction of ${\cal Z}_q$
in the simpler cases $k=1$, and $k=2$ to a group
algebra $\C[Q]$ with $Q$ being the root lattice of the
$sl(2)$ Lie algebra
\cite{FrJi88}, and
to a tensor product of  $\C[Q]$ with a quantum deformation of
a Clifford algebra \cite{Ber89}, respectively.
Then, we provide an explicit construction of the basis of the
simplest
${\cal Z}_q$ modules, and hence $\uq$ modules, for arbitrary non-zero
level $k$. These are the so called generalized Verma modules
\cite{Lep79,LePr84}. We show that the spanning vectors of
the basis  of a  generalized Verma module
as constructed in the `classical' ${\cal Z}$ algebra case in
\cite{LePr84,LePr85}, do
still form a basis for a ${\cal Z}_q$-generalized Verma module.
In section 5, we
give more quantum generalized
commutation relations satisfied by polynomials of  ${\cal Z}_q$ elements.
We think
that these relations will eventually be useful in the explicit
constructions
of all the ${\cal Z}_q$-modules, and from them
all the $\uq$-modules, and especially the standard ones
with arbitrary non-zero level $k$. Finally, Section 6
is devoted to our conclusions.

\section{The $\uq$ quantum affine algebra}

The $\uq$ affine algebra  is  a unital associative algebra
with elements $\{e_{\pm\alpha_i},  k_i^{\pm},
q^{\pm d}; \\i=0,1\}$ and defining relations in the homogeneous
gradation \cite{Dri85,Jim85}
\bac
&&{[k_i,k_j]} =0,\quad k_ik^{-1}_i=k^{-1}_{i}k_{i}=1,\\
&&q^{-d}q^d=q^dq^{-d}=1,\\
&&q^dk_iq^{-d}=k_i,\quad q^d
e_{\pm\alpha_i}q^{-d}=q^{\pm\delta_
{i0}}e_{\pm\alpha_i},\\
&& k_ie_{\pm\alpha_j} k^{-1}_i =q^{\pm (\alpha_i,
\alpha_j)}e_{\pm\alpha_i},\\
&&{[e_{\alpha_i},e_{\alpha_j}]} =
\delta_{ij}{k_i-k_i^{-1}\over q-q^{-1}},\\
&&(e_{\pm\alpha_i})^3e_{\pm\alpha_j}-[3]
(e_{\pm\alpha_i})^2e_{\pm\alpha_j}e_{\pm\alpha_i}+
[3]
e_{\pm\alpha_i}e_{\pm\alpha_j}(e_{\pm\alpha_i})^2-
e_{\pm\alpha_j}(e_{\pm\alpha_i})^3=0,
\label{qafin}
\ea
where $[x]=(q^{x}-q^{-x})/(\q)$,
$q=e^{t/2}$ is  a complex number called
a deformation parameter, and
$\{\alpha_i, i=0,1\}$ is the set of positive simple roots
of $\LS$ affine Lie algebra with the invariant symmetric
bilinear form $(\alpha_i,\alpha_j)=a_{i,j}$
(see Section 4).
Here  $a_{i,i}=2$ and $a_{i,1-i}=-2$, $i=0,1$ are the elements of
the
$\LS$ affine Cartan matrix.
Note that the
special element $\gamma=k_0k_1$
is in the center of $\uq$ and acts as $q^k$ on its highest weight
representations, with $k$ referred to as the level.
We also refer
to the above set of elements generating $\uq$ as the Chevalley basis.

The Chevalley basis consists of elements associated with
the simple roots only. One would like to describe
the commutation relations of
all basis elements associated with the infinite-dimensional
set of roots $\{\pm \alpha+n\delta; n\in \Z\}\cup\{n\delta;
n\in \Z\backslash\{0\}\}$, with $\al=\al_1$ and $\delta=
\al_0+\al_1$, and where the Serre relation
(i.e., the last relation in (\ref{qafin})) becomes redundant.
Drinfeld succeeded in finding such a basis,
which we refer as the Drinfeld basis \cite{Dri86}.
It is generated by the      elements $\{x_{n}^{\pm},
\al_m, K^{\pm}, q^{\pm d}, \gamma^{\pm 1/2}; n\in \Z,
m\in\Z^*=\Z\backslash \{0\}\}$
with defining relations
\bea
&&\g^{1/2}\g^{-1/2}=\g^{-1/2}\g^{1/2}=1,\quad
{[\g^{\pm 1/2}, y]}=0,\quad \forall y\in \uq,\label{Eq1}\\
                       &&KK^{-1}=K^{-1}K=1,\label{Eq2}\\
&&K\al_nK^{-1}=  \al_n,\label{Eq3}\\
&& Kx^\pm_n K^{-1} =q^{\pm 2} x^\pm_n,\label{Eq4}\\
&&q^{d}q^{-d}=q^{-d}q^{d}=1,\quad Kq^{\pm d}K^{-1}=q^{\pm d},
\label{Eq5}\\   &&q^dx^\pm_nq^{-d}=q^nx^\pm_n,\label{Eq6}\\
&&q^d\al_nq^{-d}=q^n\al_n,\label{Eq7}\\
&&{[\al_n,\al_m]} = {(q^{2n}-q^{-2n})(\g^{n}-\g^{-n}) \ov
nt^{2}}\delta_{n+m,0},\label{Eq8}\\
&& {[\al_n,x^{\pm}_m]}=
\pm\sqrt{2}{\g^{\mp |n|/2}(q^{2n}-q^{-2n})\over 2nt}
x^\pm_{n+m},\label{Eq9}\\
&&{[x^+_n,x^-_m]} = {\g^{(n-m)/2}\Psi_{n+m}-\g^{(m-n)/2}
\Phi_{n+m}\over q-q^{-1}},\label{Eq10}\\
&&x^\pm_{n+1}x^\pm_m-q^{\pm 2}x^\pm_mx^\pm_{n+1}=
 q^{\pm 2}x^\pm
_nx^\pm_{m+1}-x^\pm_{m+1}x^\pm_n, \label{Eq11}
\eea
where
$\Psi_n$ and $\Phi_n$
are given by
the mode expansions of the fields
$\Psi(z)$ and $\Phi(z)$, which are themselves defined by
\bac
\Psi(z)&=&\sum\limits_{n\geq 0}\Psi_nz^{-n}=K
\exp\{t\sum\limits_{n>0}\al_nz^{-n}\},\\
\Phi(z)&=&\sum\limits_{n\leq 0}\Phi_nz^{-n}=K^{-1}
\exp\{-t\sum\limits_{n<0}\al_nz^{-n}\}.
\label{algebra}
\ea
Here $z$ is a formal variable and
\be
K=\Psi_0=\Phi_0^{-1}\equiv q^{\al_0},
\label{KP}\ee
where we mean identification by the symbol $\equiv$.

The isomorphism $\rho$ between $\uq$ in the Chevalley basis
and $\uq$ in the Drinfeld basis is given explicitly by
\bac
\rho:\quad k_0\ra \g K^{-1},\\
\rho:\quad k_1\ra K,\\
\rho:\quad e_{\pm\alpha_1}\ra x^{\pm}_0,\\
\rho:\quad e_{\alpha_0}\ra x^-_{1}K^{-1},\\
\rho:\quad e_{-\alpha_0}\ra Kx^+_{-1}.
\ea

For later purposes we will use the formal variable  approach
\cite{Freal88} (instead of the usual operator product
expansion method) to    re-express
the algebra as a quantum current algebra with elements
$\{\Psi(z), \Phi(z), x^{\pm}(z), \g^{\pm 1/2},
q^{\pm d}\}$, where \cite{FrJi88}
\be
x^\pm(z)=\sum_{n\in\Z}x^\pm_nz^{-n},
\ee
and with defining relations
\bea
\g^{1/2}\g^{-1/2}=\g^{-1/2}\g^{1/2}&=&1,\quad {[\g^{\pm 1/2}, y]}=0,
\quad \forall y\in \uq,\label{eq1}\\
{[\Psi(z),\Psi(w)]}&=&0,\label{eq2}\\
{[\Phi(z),\Phi(w)]}&=&0,\label{eq3}\\
\Psi(z)\Phi(w)&=&g(wz^{-1}\g)g(wz^{-1}\g^{-1})^{-1}
\Phi(w)\Psi(z), \label{eq4}\\
\Psi(z)x^{\ep}(w)&=&
g(wz^{-1}\g^{-\ep /2})^{-\ep}x^\ep(w)\Psi(z),\label{eq5} \\
\Phi(z)x^{\ep}(w)&=&g(zw^{-1}\g^{-\ep/2})^{\ep}
x^\ep(w)\Phi(z),\label{eq6}\\
{[x^\ep(z),x^{-\ep}(w)]}&=&\ep{\delta(zw^{-1}\g^{-\ep})
\Psi(w\g^{\ep/2})- \delta(zw^{-1}\g^\ep)\Phi(z\g^{\ep/2})\over
\q},\label{eq7}\\
(z-w q^{2\ep})x^{\ep}(z)x^{\ep}(w)&=&(z q^{2\ep}-w)
x^{\ep}(w)x^{\ep}(z),\label{eq8}\\
q^{d}x^{\ep}(z)&=&x^\ep(zq^{-1})q^d,\label{eq9}\\
q^{d}\Psi(z)&=&\Psi(zq^{-1})q^d,\label{eq10}\\
q^{d}\Phi(z)&=&\Phi(zq^{-1})q^d.
\label{eq11}
\eea
Here $\ep=\pm 1$ and $g(z)$ is meant to be the following
formal power series      in $z$:
\be
g(z)=\sum_{n\in \Z_+}c_nz^n,
\ee
where the coefficients $c_n, n\in \Z_+$ are determined from the Taylor
expansion of the function
\be
f(\xi)={q^2\xi-1\over \xi-q^2}=\sum_{n\in \Z_+}c_nz^n
\ee
at $\xi=0$. Note that for this reason $f(\xi^{-1})=f(\xi)^{-1}$
but $g(z^{-1})\neq g(z)^{-1}$. In fact, $g(z)^{-1}$ is
obtained from $f(\xi)^{-1}$ in the same manner as $g(z)$ is obtained
from $f(\xi)$ \cite{FrJi88}.  In the above relations we have
also introduced the famous
$\delta$-function $\delta(z)$ which is defined as the
formal Laurent series
\be
\delta(z)=\sum_{n\in \Z}z^n,
\ee
and which plays a key role in the formal calculus approach
\cite{Freal88}. Its main properties are summarized by the
following relations:
\bac
&&\delta(z)=\delta(z^{-1}),\\
&&\delta(z)={1\over z}+{z^{-1}\over 1-z^{-1}},\\
&&G(z,w)\delta(azw^{-1})=G(z,az)\delta(azw^{-1})=
G(a^{-1}w,w)\delta(azw^{-1} ),\quad a\in \C^*,
\label{del}\ea
where $G(z,w)$ is any operator with a
 formal  Laurent  expansion in $z$ and $w$ given by
\be
G(z,w)=\sum_{n,m\in \Z}G_{n,m}z^nw^m.
\ee
Note that it is crucial that both $\delta(z)$ and
$G(z,w)$ have expansions in integral powers of $z$ and
$w$, otherwise the above properties of the $\delta$-function will not  hold.

The three relations  (\ref{eq9}), (\ref{eq10}) and (\ref{eq11})
translate the fact
that $x^{\pm}_n$, $\Psi_n$ and $\Phi_n$ are
homogeneous of the same degree $n$.

\section{The ${\cal Z}_q$
algebra}

It is well known that the Heisenberg subalgebra of $\LS$
plays a crucial role in the construction
of vertex operators and highest weight representations.
One would like to extend this role to the quantum case.
Note that in the sequel, the unit element  and
$q^{\pm d}$ are meant to be in all
the (sub)algebras defined below, so we will not consider
them unless stated otherwise.
Let $U_q(\hat{h})$ be the quantum analogue  of the
enveloping Heisenberg algebra, referred to as $q$-Heisenberg algebra.
It is a subalgebra of $\uq$ generated by $\{\alpha_n,
\g^{\pm 1/2}, n\in \Z^* \}$ with relations:
\bea
{[\al_n, \g^{\pm 1/2}]}&=&0,\quad \g^{1/2}\g^{-1/2}=\g^{-1/2}\g^{1/2}=1,\\
{[\al_n,\al_m]}& =& \delta_{n+m,0}{(q^{2n}-q^{-2n})(\g^{n}-
\g^{-n})\over 2nt^{2}}.
\label{qH}\eea
Let $U_{q}(\hat{h^{+}})$ and $U_{q}(\hat{h^{-}})$ denote the commutative
subalgebras of $U_{q}(\hat{h})$ generated by
$\{\al_n, \g^{\pm 1/2};\\  n>0\}$ and $\{\al_n; n<0\}$ respectively. By the
quantum analogue  of Poincar\'e-Birkhoff-Witt theorem  for
$U_{q}(\hat{h})$, we have
\be
U_{q}(\hat{h})=U_{q}(\hat{h^{+}})U_{q}(\hat{h^{-}}),
\ee
and consequently, the following induced $\hat{h}$-module:
\be
I(q^k)=U_q(\hat{h})\otimes_{U_q(\hat{h^{+}})}
\C[q^k]
\ee
is irreducible and isomorphic to $U_q(\hat{h^{-}})$ and
hence to the symmetric algebra $S(\hat{h^{-}})$
\cite{FrJi88}. In this
formula, $\C[q^k]$ denotes the field of complex numbers
considered as the one-dimensional $U_q(\hat{h^{+}})$-module
and on which $\g$ acts as multiplication by
$q^{k}$ and  $\al_n, n>0$ acts trivially. This means that
$S(\hat{h^{-}})$ is a  canonical $U_q(\hat{h})$-module on
which $\g$ acts as multiplication by $q^{k}$,
$\al(n)$ $(n<0)$ acts as a creation (multiplication) operator,
and $\al(n)$ ($n>0$) acts as an (derivation) annihilation operator
satisfying the relation (\ref{qH}). Henceforth, $\al(n)$
denotes the generator $\al_n$ on $S(\hat{h^{-}})$.
More precisely, the latter actions
are given by:
\bac
&&\gamma^{\pm 1}:\quad x\ra q^{\pm k}x,\\
&&\al(n):\quad x\ra \al_n x,\quad n<0,\\
&&\al(n):\quad x\ra [\al_n,x],\quad n>0,
\ea
where $x$ is any element in $\Sh$.  Moreover, the action of
$q^{\pm d}$ on $\Sh$ is defined by
\be
q^{\pm d}:\quad x\ra q^{\pm d}xq^{\pm d}.
\ee
We also denote by $\Psi(n)$ and $\Phi(n)$ the generators
$\Psi_n$ and $\Phi_n$ on $S(\hat{h^{-}})$ respectively,
that is,  $\Psi(n)$ and $\Phi(n)$ are related to
$\al(n)$ in the same way as $\Psi_n$ and $\Phi_n$ are related
to $\al_n$ (see (\ref{algebra}). The action of $\al(0)$ and hence the
actions of $\Psi(0)=q^{\alpha(0)}$ and
$\Phi(0)=q^{-\alpha(0)}$ on  $S(\hat{h^{-}})$ will be defined later.

Now we would like  to show that  the
highest weight modules of the whole quantum affine algebra
$\uq$ must be constructed as  tensor products of the form
$S(\hat{h^{-}})\otimes W$.
Here $W$ are certain vector spaces to be defined later and
which are trivial as $U_q(\hat{h})$-module. {\it
La raison d'\^etre}
of $W$ stems from the fact that  $S(\hat{h^{-}})$ is only
a $U_q(\hat{h})$-module and in general cannot be upgraded to
a $\uq$-module, which is especially true here since we are
considering $\uq$ in the homogeneous gradation. This is because
it is well
known that the $\LS$-highest weight modules in the homogeneous
gradation are constructed
as the above tensor products with non-trivial $W$ spaces
\cite{LePr84}, and so
we expect this to be also true for the $\uq$-highest weight modules.
Therefore, we have to ``correct"  $S(\hat{h^{-}})$ by tensoring
it  with additional new spaces which do not overlap with it,
so that the  resulting tensor product remains as a
$U_q(\hat{h})$-module.      Of course, this
correction will be performed by considering a minimum number
of extra spaces.

It is well known in the case of affine algebras that these
constructions can be achieved by means of  vertex
operators. The most famous vertex construction of the quantum
affine algebras is the Frenkel-Jing one, which is however
valid only for
the simply laced algebras and with the central element
$\g$ acting as $q$ (i.e., $k=1$) on their highest
weight modules \cite{FrJi88}. In this case, which will be recovered
explicitly
later for $\uq$ when we set $k=1$ in our general construction,
it turns out that the required
minimal correction
consists in tensoring the symmetric algebra with a group
algebra associated with the root lattice of the
Lie algebra corresponding to the quantum
affine algebra in question.  However,
if $k>1$  the latter correction is not
sufficient in the sense that it requires  extra spaces.
For example, if $k=2$ one needs to consider, in addition to the group
 and symmetric algebras, an exterior (Clifford) algebra. Such a
construction has been achieved in the case
of $U_q(so(2n+1))$ with level 1 by Bernard \cite{Ber89}. For
$k>2$ it was shown in \cite{BoLu94} that one needs to introduce
besides the group and symmetric algebras, a certain quantum
parafermionic algebra (though the representation theory was not
discussed there). In fact, when $k>2$ we find it easier for our
purposes here to discuss
the representation theory of the whole ${\cal Z}_q$ algebra rather than
its decomposition as a tensor product of a group algebra and
a parafermionic algebra.

Although we are concerned with $\uq$ for  $k>1$,
 the form of the vertex operators used by Frenkel and Jing
for $k=1$  led us to introduce the following vertex operators
\be
S_{\ep}^{\pm}(z)=\exp\{\pm \ep t\sum_{n>0}{\al(\pm n)\over
q^{nk}-q^{-nk}}q^{-\ep nk/2}z^{\mp n}\},\quad\ep=\pm,
\label{S}\ee
which are viewed as  formal Laurent series in $z$ with coefficients
acting on  $S(\hat{h^{-}})$. Strictly speaking, when $k=1$ the above
vertex operators reduce to the inverse operators of those
considered by Frenkel and Jing. Using  (\ref{qH}) and the
usual formal rule
\be
e^{A}e^{B}=e^{B}e^{A}e^{[A,B]}\quad{\rm if}\quad [A,[A,B]]=
[B,[A,B]]=0,
\label{rule1}\ee
for some operators $A$ and $B$, we find
\bac
S_{\ep}^{+}(z)S_{\epp}^{-}(w)&=&\displaystyle
{(q^{k-2-(\ep+\epp)k/2}wz^{-1};q^{2k})_{\infty}^{\ep\epp}\over
(q^{k+2-(\ep+\epp)k/2}wz^{-1};q^{2k})_{\infty}^{\ep\epp}}
S_{\epp}^{-}(w)S_{\ep}^{+}(z),\\
S_{\ep}^{\pm}(z)S_{\epp}^{\pm}(w)&=&
S_{\epp}^{\pm}(w)S_{\ep}^{\pm}(z),
\label{SS}\ea
where as usual $(x;y)_{\infty}$ means
\be
(x;y)_{\infty}=\prod_{n=0}^{\infty}(1-xy^{n}).
\ee
Each factor $(1-wz^{-1}q^x)^{-1}$ in (\ref{SS}) is understood
as the formal power series $\sum_{n\geq 0}w^{n}z^{-n}q^{xn}$.
Moreover, using (\ref{qH}) and the formal rule
\be
{[A, e^{B}]}={[A,B]}e^{B},\quad {\rm if}\quad [B,[A,B]]=0,        \ee
for some operators $A$ and $B$, we obtain
\bac
{[\al(n),S^+_{\ep}(z)]}&=&0,\\
{[\al(-n),S^+_{\ep}(z)]}&=&\displaystyle
-{\ep q^{-\ep nk/2}z^{-n}(q^{2n}-      q^{-2n})\over nt}S^+_{\ep}(z),\\
{[\al(-n),S^-_{\ep}(z)]}&=&0,\\
{[\al(n),S^-_{\ep}(z)]}&=&\displaystyle
-{\ep q^{-\ep nk/2}z^{n}(q^{2n}-q^{-2n})\over nt}S^-_{\ep}(z),
\label{aS}
\ea
where $n>0$ and $\ep=\pm$. For future purposes let us note here
that one can easily show that the commutation relations (\ref{eq5})
and (\ref{eq6}) are  equivalent to
\be
{[\al_n,x^{\ep}(z)]}={\ep q^{-\ep |n|k/2}z^{n}(q^{2n}-
q^{-2n})\over nt}x^{\ep}(z),\quad \ep=\pm;\quad n\in \Z\backslash  \{0\}.
\label{ax}\ee
Let us now define the main new objects of this paper,
which we refer to as the ``${\cal Z}_q$ operators ${\cal Z}^{\ep}_n$,"
as the Laurent modes in     \be
 {\cal Z}^{\ep}(z)=\sum_{n\in \Z}{\cal Z}_{n}^{\ep}z^{-n},
\ee
where
\be
{\cal Z}^{\ep}(z)=S^-_{\ep}(z)x^{\ep}(z)S^+_{\ep}(z).
\label{ZOP}
\ee
These operators ${\cal Z}^{\ep}_n$ are the quantum analogues of
the (classical, $q=1$) ${\cal Z}$-operators that have been extensively
studied in the literature (see for example Reference
\cite{LePr85}).
By abuse of terminology, we refer  also to the ``currents
${\cal Z}^{\ep}(z)$"
as ${\cal Z}_q$ operators,  but strictly speaking they are the generating
functions of the latter operators.
Let us stress here that $z$ is a formal variable, $\Z$ is the
set  of integer numbers, whereas
${\cal Z}^{\ep}_n$ are operators so there should not be any
confusion with these notations.
Let us also denote by ${\cal Z}_q$ the algebra
generated by $\{\Psi_0,\Phi_0,{\cal Z}^\ep_n;n\in \Z\}$.
The defining relations of this
algebra, which we
refer to as ``the quantum generalized commutation relations"
will be given shortly below.

The space $W$ on which this algebra acts non-trivially is the necessary
space to be tensored with $S(\hat{h^{-}})$
such that $S(\hat{h^{-}})\otimes W$ is a $\uq$-module.
Therefore  $W$ will be defined if we know all the properties of
the ${\cal Z}_q$ operators, that is, their  relations
with $\uq$ itself and their algebra (i.e., quantum generalized
commutation relations).

For this purpose, let ${\cal Z}(\ep|n)$ and $X(\ep|n)$  represent
the actions  of
${\cal Z}^{\ep}_n$ and $x^{\ep}_n$ on $W$ and $\Sh\otimes W$
respectively. By definition, the relation between the latter operators
is given through their generating functions by
(\ref{ZOP}). Next, it can easily be checked that the relations (\ref{aS})
and  (\ref{ax}) imply
that the ${\cal Z}_q$ operators commute with the quantum
Heisenberg algebra $U_q(\hat{h})$, i.e.,
\bac
{[\al(n),{\cal Z}(\ep|z)]}&=&0,\quad n\in \Z\backslash\{0\},\\
{[\g^{\pm 1/2},{\cal Z}(\ep|z)]}&=&0,\quad \ep=\pm.
\label{aZ}
\ea
This is a very important result, which is in fact
the main motivation behind the particular choice for the
forms of  $S^{\ep}_{\epp}(z)$ as given by (\ref{S}).
This is because the
symmetric algebra $S(\hat{h^{-}})$ realizes
already the quantum Heisenberg subalgebra $U_q(\hat{h})$ and since
the ${\cal Z}_q$ operators commute with $U_q(\hat{h})$,
we can define then the actions of $U(\hat{h})$ and
${\cal Z}_q$ as follows on $S(\hat{h^{-}})\otimes W$:
\bac
x: u\otimes v\ra xu\otimes v,\\
y: u\otimes v\ra u\otimes yv,
\ea
where $x\in U_q(\hat{h})$, $y\in {\cal Z}_q$,
$u\in S(\hat{h^{-}})$,   and $v\in W$. From the
Laurent expansion in $zw^{-1}$ of both sides of (\ref{eq5})
and (\ref{eq6}), we obtain  the relations
\bac
\Psi(0)X(\ep|w)&=&
q^{2\ep}X(\ep|w)\Psi(0),\\
\Phi(0)X(\ep|w)&=&q^{-2\ep}X(\ep|w)\Phi(0),
\label{0PX}\ea
which, because of the relation (\ref{Eq4}), amount to
\bac
\Psi(0){\cal Z}(\ep|w)&=&
q^{2\ep}{\cal Z}(\ep|w)\Psi(0),\\
\Phi(0){\cal Z}(\ep|w)&=&q^{-2\ep}{\cal Z}(\ep|w)\Phi(0).
\label{PZ}\ea
Combining these relations with (\ref{aZ}) and
\bac
\Psi(z)&=&\Psi_{0}
S^{+}_{\ep}(zq^{-3\ep k/2})S^{+}_{-\ep}(zq^{3\ep k/2}),\quad \ep=  \pm,\\
\Phi(z)&=&
\Phi_{0}S^{-}_{\ep}(zq^{3\ep k/2})S^{-}_{-\ep}
(zq^{-3\ep k/2}),\quad \ep=\pm,
\label{ps}\ea
we arrive finally at
\bac
\Psi(n){\cal Z}(\ep|w)&=&
q^{2\ep}{\cal Z}(\ep|w)\Psi(n),\quad n\geq 0, \\
\Phi(n){\cal Z}(\ep|w)&=&q^{-2\ep}{\cal Z}(\ep|w)\Phi(n),\quad n\leq 0.
\label{ppz}\ea
In  the sequel, however, we will only consider the relations
(\ref{aZ}) and (\ref{PZ}) but not (\ref{ppz}) since the latter
is an immediate consequence of the former and (\ref{ps}).

Clearly, the action of the generating functions
\be
X(\ep|z)=\sum_{n\in \Z}X(\ep|n)z^{-n}
\ee
on $S(\hat{h^{-}})\otimes W\otimes \C[z,z^{-1}]$ decomposes then as  \be
X(\ep|z)=S^-_{-\ep}(zq^{-\ep k})S^+_{-\ep }(zq^{\ep k})
\otimes {\cal Z}(\ep|z),\quad \ep=\pm,
\label{XS}\ee
where we have used
\be
\bigl(S^{\pm}_{\ep}(z)\bigr)^{-1}=
S^{\pm}_{-\ep}(zq^{\pm\ep k}),\quad \ep=\pm,
\label{PSS}
\ee
and (\ref{aS}) to express  $X_{\ep}(z)$ in terms of ${\cal Z}(\ep|z)$.
We now define the actions of $q^{\pm d}$ and
$\gamma^{\pm 1}$ on $S(\hat{h^{-}})\otimes W$ as
\bac
q^{\pm d}: u\otimes v\ra q^{\pm d}u\otimes q^{\pm d}v,\\
\g^{\pm 1}: u\otimes v\ra q^{\pm k}u\otimes v,
\ea
where $u\in \Sh$ and $v\in W$.
The relation between ${\cal Z}(\ep|z)$ and $q^{\pm d}$, which reads
as
\be
q^{d}{\cal Z}(\ep|z)={\cal Z}(\ep|zq^{-1})q^d,
\label{hab}\ee
 can  easily be
derived from (\ref{eq9}) and
\be
q^{d}S^{\ep}_{\epp}(z)=S^{\ep}_{\epp}(zq^{-1})q^{d},
\ee
which, in turn, can be obtained from (\ref{rule1}).
Relation (\ref{hab}) means that the  ${\cal Z}_q$ algebra is graded
(it inherits the
gradation of $\uq$) and that the ${\cal Z}_q$ operators ${\cal Z}(\ep|n)$ are
homogeneous of degree $n$.

Let us now turn to the derivation of  the defining
relations (besides (\ref{KP}) and (\ref{PZ})) of the ${\cal Z}_q$ algebra,
that is, the quantum
generalized commutation relations. They simply follow from
the substitution of $X(\ep|z)$ as given by (\ref{XS}) in
(\ref{eq7}) and (\ref{eq8}). We find
\bea
&&{(q^{k+2}wz^{-1};q^{2k})_{\infty}\over
(q^{k-2}wz^{-1};q^{2k})_{\infty}}{\cal Z}(\ep|z){\cal Z}(-\ep|w)-
{(q^{k+2}zw^{-1};q^{2k})_{\infty}\over
(q^{k-2}zw^{-1};q^{2k})_{\infty}}{\cal Z}(-\ep|w){\cal Z}(\ep|z)\nonumber\\
&&={\ep\over q-q^{-1}}
\left(\Psi_{0}\delta(zw^{-1}q^{-\ep k})
-\Phi_{0}\delta(zw^{-1}q^{\ep k})\right)\nonumber\\
&&={1\over q-q^{-1}}
\left(q^{\ep \al(0)}\delta(zw^{-1}q^{- k})
-q^{-\ep\al(0)}\delta(zw^{-1}q^{ k})\right),\\
\label{qZ1}
&&(z-q^{2\ep}w){(q^{k-2-k\ep}wz^{-1};q^{2k})_{\infty}\over
(q^{k+2-k\ep}wz^{-1};q^{2k})_{\infty}}{\cal Z}(\ep|z){\cal Z}(\ep|w)\nonumber\\
&&=(q^{2\ep}z-w){(q^{k-2-k\ep}zw^{-1};q^{2k})_{\infty}\over
 (q^{k+2-k\ep}zw^{-1};q^{2k})_{\infty}}{\cal Z}(\ep|w){\cal Z}(\ep|z).
\label{qZ2}
\eea
In summary, in addition to the latter quantum generalized commutation
relations, the
${\cal Z}_q$ fields ${\cal Z}(\ep|z)$ must satisfy the following relations
with
the elements of $\uq$ when they act on the $W$ part of the
tensor product $S(\hat{h^{-}})\otimes W$:
\bea
{[\al(n),{\cal Z}(\ep|z)]}&=&0,\quad n\in \Z\backslash\{0\},
\label{R1}\\
\Psi(0){\cal Z}(\ep|z)&=&
q^{2\ep}{\cal Z}(\ep|z)\Psi(0),
\label{R2}\\
\Phi(0){\cal Z}(\ep|z)&=&q^{-2\ep}{\cal Z}(\ep|z)\Phi(0),\label{R3}\\
X(\ep|z)&=&S^-_{-\ep}(zq^{-\ep k})S^+_{-\ep }
(zq^{\ep k})\otimes {\cal Z}(\ep|z),\label{R4}\\
q^{d}{\cal Z}(\ep|z)&=&{\cal Z}(\ep|zq^{-1})q^d,\label{R5}\\
{[\g^{\pm},{\cal Z}(\ep|z)]}&=&0.
\label{R6}\eea
All these relations will be useful in the explicit construction
of the space $W$ from the operators ${\cal Z}(\ep|n)$.
This is illustrated in the next section.
\section{Explicit constructions of some $\uq$-modules}

Let us briefly recall the definition of some
$\uq$-modules \cite{Idzu93,JiMi94,LePr84}.
For this, we still
need some notions from $\LS$ affine
algebra, which is generated by
$\{e_i,f_i,h_i,d; i=0,1\}$. We define on its Cartan
subalgebra $\hat{h}=\C h_0+\C h_1+\C d$ an invariant symmetric bilinear
form $(\quad,\quad)$ by
\be
(h_i,h_j)=a_{i,j},\quad (h_i,d)=\delta_{i,0},\quad
(d,d)=0,
\ee
where
\be
a_{i,i}=2,\quad a_{i,1-i}=-2,\quad i=0,1,
\ee
are the elements of the $\LS$ Cartan matrix. Let
$\hat{h}^*=\C \Lambda_0+\C \Lambda_1+\C \delta=
\C \al_0+\C \al_1+\C \Lambda_0$ be the dual space to
$h$ with
\be
<\Lambda_i,h_j>=\delta_{i,j},\quad
<\delta,d>=1,\quad
<\Lambda_i,d>=0,\quad
<\delta, h_i>=0,
\ee
where
\be
<\quad,\quad>:\quad \hat{h}^*\otimes \hat{h}\ra \C.
\ee
is the natural pairing, the vectors $\Lambda_i$ are the
fundamental weights, the vectors $\al_i$ are the positive
roots and $\delta=\al_0+\al_1$ is the null root. One can induce
a symmetric bilinear form $(\quad,\quad)$ on $\hat{h}^*$ by
\be
(f(x),f(y))=(x,y),\quad f(x),f(y)\in
\hat{h}^*\quad x,y\in \hat{h}.
\ee
Here for a given $x\in \hat{h}$, $f(x)$ is defined to be
the unique vector in $\hat{h}^*$ such that
\be
<f(x),y>=(x,y),\quad \forall y\in \hat{h}.
\ee
This relation allows the identification of any element in
$\hat{h}$ with a unique element in $\hat{h}^*$. In particular,
the elements $h_i$ and $d$ of the Cartan subalgebra are
identified with $\al_i$ and $\Lambda_0$ respectively.
Note that this identification should note be confused with
the identification $h_1=\al(0)$, which is the classical
analogue of the identification that we have made implicitly
in $K=q^{h_1}=q^{\al_1(0)}=q^{\al(0)}$ in the case of $\uq$ (see
(\ref{KP}). To see the
relation between the two identifications, let $v_{\lambda}$
be a weight vector of $\hat{h}$ with weight $\lambda$ then
\be
h_1 v_{\lambda}= <\lambda,h_1>v_\lambda\equiv
<\lambda,\al(0)>v_\lambda\equiv (\lambda,\al_1)v_\lambda\equiv
(\lambda,\al)v_\lambda.
\ee

More explicitly the symmetric bilinear form on $\hat{h}^*$
is given  by
\bac
&&(\Lambda_i,\Lambda_j)={1\over 2}\delta_{i,1}\delta_{j,1},
\quad  (\Lambda_i,\delta)=1, \quad (\delta,\delta)=0,\\
&&(\al_i,\al_j)=a_{i,j},\quad (\al_i,\Lambda_0)=
\delta_{i,0},\quad (\Lambda_0,\Lambda_0)=0,\quad i,j=0,1.
\ea

The weights $\lambda\in \hat{h}^*$ such that
\be
\lambda=n_{0}\Lambda_{0}+n_{1}\Lambda_{1},\quad n_{0},n_{1}\in \N,
\ee
are called dominant integral weights, and $n_0+n_1=k$
is  the level that we have introduced previously.

As defined in Section 2, the algebra
$\uq$ is generated by $\{e_i,f_i,K^{\pm 1},\gamma^{\pm 1},
q^{\pm d};\\
i=0,1\}$.
Let $V$ be a $\uq$-module and $\mu\in \hat{h}$, the subspace
$V_\mu\subset V$ defined by
\be
V_\mu =\{v\in V/ K^{\pm 1}v=q^{\pm <\mu,h_1>}v,\quad\gamma^{\pm 1}
v=q^{\pm k}v,\quad q^{\pm d}v=q^{\pm <\mu ,d>}v\},
\ee
is called a $\mu$-weight space, and any $v\in V_\mu$ is referred to
as a $\mu$-weight vector.
The module $V$ becomes a weight
module if it  is the direct sum of its weight spaces.
A $\uq$-highest weight vector  $v_\lambda$ in $V$ is a
$\lambda$-weight vector which satisfies the
additional condition
\be
e_iv_\lambda=0,\quad i=0,1.
\ee
The space $V$ is called a
$\uq$-highest weight module if it generated from a
$\lambda$-highest weight vector $v_\lambda$. In this case,
$V$ is also a $\uq$ weight module and $v_\lambda$ is unique
(up to a multiplication by a scalar), and hence we label $V$ by the
weight $\lambda$ as $V(\lambda)$.
The $\uq$-module $V(\lambda)$ is called standard if it is
generated from a highest weight vector $v_\lambda$
with a dominant integral weight  $\lambda$ and such that
\be
f_i^{<\lambda,h_i>+1}v_\lambda=0,\quad i=0,1;
\ee
in which case it is irreducible, and thus  called also an
irreducible highest weight module.
Let $U_q(sl(2))$ be a subalgebra of $\uq$
generated by $\{e_1,f_1,K^{\pm 1}\}$ and
let $M=\C v_0$ be the trivial one-dimensional
$U_q(sl(2))$-weight module.
We also introduce $\uq_{\geq 0}$, $\uq_{>0}$ and $\uq_{<0}$ as
three subalgebras of $\uq$
generated by all elements with nonnegative, positive and
negative degrees with respect to $q^{d}$, respectively.
We equip $M$ with a $\uq_{\geq 0}$-module structure by
\bac
q^{\pm d}v_0&=&q^{\pm a}v_0,\quad a\in\C,\\
\gamma^{\pm 1}v_0&=&q^{\pm k}v_0,\\
xv_0&=&0,\quad \forall x\in \uq_{>0}.
\ea
Here $a$ is scalar, which can be set to 0, without loss of
generality. From $M$, we induce the following
$\uq$-module:
\be
G(M)=\uq\otimes_{\uq_{\geq 0}}M,
\ee
which is called a generalized Verma module \cite{Lep79,LePr84}.
In fact, since $M$ is a one-dimensional $U_q(sl(2))$-module,
$G(M)$ is the simplest example of a $\uq$-generalized Verma module.
Slightly more complicated example of generalized
Verma modules can be constructed from higher dimensional
$M$ modules. We define the module $W(M)$
through the following isomorphism of
$\uq$-modules:
\be
G(M)\simeq \uq_{<0}\otimes_{\C}M\simeq \Sh\otimes_{\C} W(M).
\ee
Both modules $G(M)$ and $W(M)$ are $q^d$-weight modules with
weight space decompositions:
\bac
G(M)&=&\oplus_{n\leq 0}G(M)_n,\\
W(M)&=&\oplus_{n\leq 0}W(M)_n.
\ea
{}From the work of Lusztig \cite{Lus88},
we know that for generic $q$ the
dimensions of weight spaces of $\uq$-modules are the
same as those of the weight spaces of the $\LS$-modules.
The characters $\chi(G(M))$ and $\chi(W(M))$
of $G(M)$ and $W(M)$ have been computed
in the $\LS$ case in \cite{LePr84,LeWi84}  and are given by:
\bac
\chi(G(M))&=&\displaystyle\sum_{n\geq 0}dim(G(M)_n)p^n=
{1\over \prod_{n>0}(1-p^{n})^{3}},\\
\chi(W(M))&=&\displaystyle\sum_{n\geq 0}dim(W(M)_n)p^n=
{1\over \prod_{n>0}(1-p^{n})^{2}},
\label{cha}\ea
where $p\in\C^*$ is a formal variable. The second character
will allow us to prove the linear independence of
 a particular set of vectors constructed from the
${\cal Z}_q$ operators
${\cal Z}(\ep|n)$ and which span
$W(M)$. This means that the latter  set of vectors
is a basis for $W(M)$.

Let us now address the explicit constructions of the standard
modules in the cases $k=1$ and $k=2$. As explained previously, we should
address only the  explicit
constructions of the space $W$ in terms of the ${\cal Z}_q$ operators
since the $\Sh$ part is already constructed in terms of
polynomials of $\al(n)$. Consequently, in the sequel we will mainly
concentrate on the $W$ part of the $\uq$-modules.\\

\noindent Case I:$\quad k=1$\hfill

In this case, the relations (\ref{qZ1})-(\ref{R6}) satisfied by the
${\cal Z}_q$ operators ${\cal Z}(\ep|z)$ simplify significantly and reduce to
 the following relations:
\bea
&&{{\cal Z}(\ep|z){\cal Z}(-\ep|w)\over
(1-q^{-1}wz^{-1})
(1-qwz^{-1})}-{{\cal Z}(-\ep|w){\cal Z}(\ep|z)\over
(1-q^{-1}zw^{-1})
(1-qzw^{-1})}\nonumber\\
&&={\ep\over q-q^{-1}}
\left(\Psi_{0}\delta(zw^{-1}q^{-\ep })
-\Phi_{0}\delta(zw^{-1}q^{\ep })\right)\nonumber\\
&&={1\over q-q^{-1}}
\left(q^{\ep \al(0)}\delta(zw^{-1}q^{-1 })
-q^{-\ep\al(0)}\delta(zw^{-1}q)\right),
\label{equa1}\\
&&w^{2}Z(\ep|z){\cal Z}(\ep|w)=z^{2}{\cal Z}(\ep|w){\cal Z}(\ep|z),
\label{equa2}\\
&&{[\al(n),{\cal Z}(\ep|z)]}=0,\quad n\in \Z\backslash\{0\},
\label{equa3}\\
&&\Psi(0){\cal Z}(\ep|z)=
q^{2\ep}{\cal Z}(\ep|z)\Psi(0),
\label{equa4}\\
&&\Phi(0){\cal Z}(\ep|z)=q^{-2\ep}{\cal Z}(\ep|z)\Phi(0),
\label{equa5}\\
&&X(\ep|z)=S^-_{-\ep}(zq^{-\ep k})S^+_{-\ep }(zq^{\ep k})
\otimes {\cal Z}(\ep|z),
\label{equa6}\\
&&q^{d}{\cal Z}(\ep|z)={\cal Z}(\ep|zq^{-1})q^d,
\label{equa7}\\
&&{[\g^{\pm},{\cal Z}(\ep|z)]}=0.
\label{equa8}
\eea
We now construct the appropriate space $W$ that is compatible
with the above relations. We start with  (\ref{equa4})
and (\ref{equa5}). These relations imply that ${\cal Z}(\ep|z)$ cannot
be trivial, i.e., it does not act like a scalar which can
be rescaled to 1 on $W$, and hence
$W$ cannot be a trivial one-dimensional space. If the latter
statement were not true,
the relations (\ref{0PX}) would mean that $X(\ep|z)$ acts trivially
 also on $\Sh$ because of the fact that both $\Psi_0$
and $\Phi_0$ do act like the identity on $\Sh$. But if $X(\ep|z)$ acts
trivially on $\Sh$, the relation (\ref{eq7}) implies that $\Psi(z)$ and
$\Phi(z)$ also act trivially on $\Sh$, which is impossible from
the definition itself of $\Sh$. Therefore ${\cal Z}(\ep|z)$ cannot act like
the identity on $W$. To make it act nontrivially and consistently
with $\Psi_0$ and $\Phi_0$, let us write the positive simple
root $\al_1=\al$ and introduce the $sl(2)$ root lattice
$Q=\Z \al$ and the weight lattice $P=\Z\al/2$, which decomposes as
$P=Q\cup(Q+\al/2)$. Denote by $\C[P]$ the group
algebra with the basis $\{e^{\beta}, \beta\in
P\}$ and the commutative multiplication
\be
e^{\beta}e^{\gamma}=
e^{\gamma}e^{\beta}=e^{\beta+\gamma},\quad \beta,\gamma\in [P].
\ee
Obviously, we have $\C[Q]\subset \C[P]$.
We define the action of $\al(0)\in End(\C[P])$ on
$\C[P]$ by
\be
\al(0):\quad e^\beta\ra (\al,\beta)e^\beta.
\ee
This implies that the formal Laurent series
$z^{\al(0)}\in End(\C[Q]\otimes [z,z^{-1}])$ acts on $\C[P]$ as
\be
z^{\al(0)}:\quad e^\beta\ra z^{ (\al,\beta)}e^{\beta}.
\ee
It can easily be checked that the latter two relations are
equivalent to
\bac
{[\al(0),e^\beta]}&=&(\al,\beta)e^\beta,\\
z^{\al(0)}e^\beta&=&z^{(\al,\beta)}e^\beta z^{\al(0)},
\label{ae}\ea
respectively. We define also the action of $q^{ d}$ on
$\C[P]$  (i.e., $\C[P]$ is a graded space) by
\be
q^{d}:\quad e^\beta\ra e^{-1/2(\beta,\beta)}e^\beta,\quad
\beta\in\C[P],
\ee
that is,
\be
q^{d}e^{\beta}=e^{\beta}q^{d}q^{-\beta(0)-(\beta,\beta)/2}.
\ee

{}From the relations (\ref{equa4}), (\ref{equa5}) and (\ref{ae})
we see that  as a candidate for ${\cal Z}(\ep|z)$ we can take
\be
{\cal Z}(\ep|z)=e^{\ep \al},\quad \ep=\pm,
\ee
where $\al$ is the $sl(2)$ positive simple root. However,
according to the relation (\ref{equa7}) this candidate does
not have the right degree since it does not depend on the formal
variable $z$.
This is corrected by taking instead the following candidate for
${\cal Z}(\ep|z)$:
\be
{\cal Z}(\ep|z)=e^{\ep \al}z^{\ep\al(0)+{1\over 2}(\al,\al)}=
e^{\ep \al}z^{\ep\al(0)+1},\quad \ep=\pm.
\ee
With this choice, it can  easily be verified that all the
remaining relations are satisfied. The less trivial one is
(\ref{equa1}), which can be shown as follows:
\bac
&&\displaystyle{\cal Z}(\ep|z){\cal Z}(-\ep|w)\left({1\over
(1-q^{-1}wz^{-1})
(1-qwz^{-1})}-{z^2w^{-2}\over
(1-q^{-1}zw^{-1})
(1-qzw^{-1})}\right)\\
&&\displaystyle={\cal Z}(\ep|z){\cal Z}(-\ep|w)\left({1\over
1-qwz^{-1}}\delta(zw^{-1}q^{-1})-{qzw^{-1}
\over
1-qzw^{-1}}\delta(zw^{-1}q)\right)\\
&&=\displaystyle{q{\cal Z}(\ep|z){\cal Z}
(-\ep|zq^{-1})\delta(zw^{-1}q^{-1})-
q^{-1}{\cal Z}(\ep|z){\cal Z}(-\ep|zq)
\delta(zw^{-1}q)
\over q-q^{-1}}\\
&&\displaystyle ={1\over q-q^{-1}}
\left(q^{\ep \al(0)}\delta(zw^{-1}q^{-1 })
-q^{-\ep\al(0)}\delta(zw^{-1}q)\right),
\ea
where we have used
(\ref{ae}), (\ref{del}), and
\be
{\cal Z}(\ep|z){\cal Z}(-\ep|w)=(zw^{-1})^{\ep\al(0)-1}.
\ee

As a conclusion,
$X(\ep|z)$  acts on $\Sh\otimes \C[P]$ as
$S^-_{-\ep}(zq^{-\ep })S^+_{-\ep }(zq^{\ep })\otimes
e^{\ep \al}z^{\ep\al(0)+1}$ and is single-valued. Moreover,
the subspaces $\Sh\otimes \C[Q]$ and $\Sh\otimes e^{\alpha\over
2}\C[Q]$ whose direct sum is $\Sh\otimes\C[P]$ are invariant and
 irreducible. They are in fact isomorphic to the standard (basic) modules
$V(\Lambda_0)$ and $V(\Lambda_1)$ with highest weight vectors
realized as $1\otimes 1$ and $1\otimes e^{\al\over 2}$, respectively.
\\

\noindent Case II:$\quad k=2$\hfil

In this case, the relations (\ref{qZ1})-(\ref{R6})
satisfied by the ${\cal Z}_q$ operators reduce to
\bea
&&\displaystyle {{\cal Z}(\ep|z){\cal Z}(-\ep|w)\over
1-wz^{-1}}-{{\cal Z}(-\ep|w){\cal Z}(\ep|z)\over
1-zw^{-1}}
={\ep\over q-q^{-1}}
\left(\Psi_{0}\delta(zw^{-1}q^{-2\ep })
-\Phi_{0}\delta(zw^{-1}q^{2\ep })\right)\nonumber\\
&&\displaystyle ={1\over q-q^{-1}}
\left(q^{\ep \al(0)}\delta(zw^{-1}q^{-2 })
-q^{-\ep\al(0)}\delta(zw^{-1}q^{2 })\right),
\label{equal1}\\
&&(z-wq^{2\ep})(1-wz^{-1}q^{-2\ep}){\cal Z}(\ep|z){\cal Z}(\ep|w)=
(zq^{2\ep}-w)(1-zw^{-1}q^{-2\ep}){\cal Z}(\ep|w){\cal Z}(\ep|z),
\nonumber\\
&&\label{equal2}\\
&&{[\al(n),{\cal Z}(\ep|z)]}=0,\quad n\in \Z\backslash\{0\},
\label{equal3}\\
&&\Psi(0){\cal Z}(\ep|z)=
q^{2\ep}{\cal Z}(\ep|z)\Psi(0),
\label{equal4}\\
&&\Phi(0){\cal Z}(\ep|z)=q^{-2\ep}{\cal Z}(\ep|z)\Phi(0),
\label{equal5}\\
&&X_{\ep}(z)=S^-_{-\ep}(zq^{-2\ep })S^+_{-\ep }
(zq^{2\ep })\otimes {\cal Z}(\ep|z),
\label{equal6}\\
&&q^{d}{\cal Z}(\ep|z)={\cal Z}(\ep|zq^{-1})q^d,
\label{equal7}\\
&&{[\g^{\pm},{\cal Z}(\ep|z)]}=0.
\label{equal8}
\eea

As in the case $k=1$, in order to satisfy the relations (\ref{equal4}) and
(\ref{equal5}) we  need to consider the group algebras $\C[P]$ and
$\C[Q]$. The operators $\al(0)$ and $z^{\al(0)}$ act
on them in the same manner as in the case $k=1$ but now $q^{d}$ acts on them
slightly differently. Its action is rescaled by a factor 2($=k$) as:
\be
q^{d}:\quad e^\beta\ra e^{-1/4(\beta,\beta)}e^\beta,\quad
\beta\in\C[P],
\ee
which is equivalent to
\be
q^{d}e^{\beta}=e^{\beta}q^{d}q^{-{\beta(0)\over 2}-(\beta,\beta)/4},\quad
\beta\in\C[P].
\ee
The candidate for ${\cal Z}(\ep|z)$ that is compatible with this
action of $q^{d}$, and the relations (\ref{equal4}) and
(\ref{equal5}) is given by
\be
{\cal Z}(\ep|z)=e^{\ep \al}z^{\ep{\al(0)\over 2}+{1\over 4}(\al,\al)}=
e^{\ep \al}z^{\ep{\al(0)\over 2}+{1\over 2}},\quad \ep=\pm.
\label{cand1}\ee
It is also trivial that this candidate satisfies both relations
(\ref{equal3}) and (\ref{equal8}). However, it can easily be checked that
the remaining relations (\ref{equal1}) and (\ref{equal2})
are not satisfied with this candidate. This means that the
group algebra $\C[P]$
is not big enough and therefore we must enlarge it with at least
an extra algebra, which will turn out to be the Clifford algebra.
To see this, let $\psi(\ep|z)$ be a generating function of the
elements of yet an unknown quantum deformation of the Clifford
algebra such that
\bac
{[\al(n),\psi(\ep|z)]}&=&0,\\
{[\Psi_0,\psi(\ep|z)]}&=&0,\\
{[\Phi_0,\psi(\ep|z)]}&=&0,\\
q^{d}\psi(\ep|z)&=&\psi(\ep|zq^{-1})q^{d},\\
\gamma\psi(\ep|z)&=&\psi(\ep|z)\gamma.
\ea
These equations guarantee that the following candidate for
${\cal Z}(\ep|z)$ still preserves  all the relations that are already
satisfied by our first candidate (\ref{cand1}):
\be
{\cal Z}(\ep|z)=e^{\ep \al}z^{\ep{\al(0)\over 2}+{1\over 4}(\al,\al)}
\psi(\ep|z)=
e^{\ep \al}z^{\ep{\al(0)\over 2}+{1\over 2}}\psi(\ep|z),\quad \ep=\pm.
\label{cand2}\ee
In terms of these fields $\psi(\ep|z)$, the remaining
two relations which must also be satisfied take the following form
after using (\ref{ae}):
\bea
&&\displaystyle{1\over q-q^{-1}}
\left(q^{\ep \al(0)}\delta(zw^{-1}q^{-2 })
-q^{-\ep\al(0)}\delta(zw^{-1}q^{2 })\right),
\nonumber\\
&&=(zw^{-1})^{-{1\over 2}}(zw^{-1})^{{\ep\over 2}\al(0)}
\left({\psi(\ep|z)\psi(-\ep|w)\over
1-wz^{-1}}-{zw^{-1}\psi(-\ep|w)\psi(\ep|z)\over
1-zw^{-1}}\right)\nonumber\\
&&=(zw^{-1})^{-{1\over 2}}(zw^{-1})^{{\ep\over 2}\al(0)}
{\{\psi(\ep|z),\psi(-\ep|w)\}\over
1-wz^{-1}},
\label{equala1}\\
&&z(z-wq^{2\ep})(1-wz^{-1}q^{-2\ep})\psi(\ep|z)\psi(\ep|w)
=
w(zq^{2\ep}-w)(1-zw^{-1}q^{-2\ep})\psi(\ep|w)\psi(\ep|z),\nonumber\\
&&\label{equala2}
\eea
where $\{A,B\}=AB+BA$ is the anticommutator of $A$ and $B$.
Clearly, (\ref{equala1}) is satisfied if the following
conditions on $\psi(\ep|z)$ hold:
\bea
\{\psi(\ep|z),\psi(-\ep|w)\}&=&\delta(zw^{-1}q^{-2})+
\delta(zw^{-1}
q^{2}),\quad {\rm if }\quad z^{\al(0)\over 2}\in
End(e^{\al\over 2}\C[Q]\otimes \C[z,z^{-1}]),\nonumber\\
&&\label{fir1}\\
\{\psi(\ep|z),\psi(-\ep|w)\}&=&(zw^{-1})^{1\over 2}
(q^{-1}\delta(zw^{-1}q^{-2})+q\delta(zw^{-1}
q^{2})),\quad {\rm if }\quad z^{\al(0)\over 2}\in End(\C[Q])
\otimes\C[z,z^{-1}].\nonumber\\
&&\label{fir2}
\eea
Since the right hand sides of these equations
do not depend on $\ep$ we might try to identify
$\psi(+|z)$
and $\psi(-|z)$. Relation (\ref{fir1})
means that $\psi(z)\equiv \psi(\ep|z)$ has
the Laurent expansion
\be
\psi(z)=\sum_{n\in \Z}\psi_nz^{-n},\quad
\ep=\pm.
\ee
Substituting this expansion back in (\ref{fir1}) and comparing the
coefficients of powers of $zw^{-1}$, we obtain the
following anticommutation relations for the modes $\psi_n$:
\be
\{\psi_{n},\psi_{m}\}=(q^{2n}+q^{-2n})\delta_{n+m,0},
\quad {\rm if }\quad z^{\al(0)\over 2}\in End(e^{\al\over
2}\C[Q]\otimes\C[z,z^{-1}]).
\label{co1}
\ee
This is the usual Rammond (R) sector for the modes $\psi_n;n\in \Z$.
Similarly (\ref{fir2}) enforces the following   expansion:
\be
\psi(z)=\sum_{r\in \Z+{1\over 2}}
\psi_rz^{-r},\quad
\ep=\pm.
\ee
As in the (R) sector, relation (\ref{fir2})
leads to the following anticommutation relations for the modes
$\psi_r$:
\be
\{\psi_{r},\psi_{s}\}=(q^{2r}+q^{-2r})\delta_{r+s,0},
\quad {\rm if }\quad z^{\al(0)\over 2}\in End(\C[Q]\otimes
\C[z,z^{-1}]).
\label{co2}\ee
This is the familiar Neuveu-Schwarz (NS) sector for the modes
$\psi_r$.
Note that both  relations (\ref{fir1}) and
(\ref{fir2}) imply that
\be
\psi(z)\psi(z)=\delta(q^2),
\ee
which means that unlike in the classical case, the divergence
in the product $\psi(z)\psi(z)$ at the same ``point" $z$
is regularized by the deformation parameter,
which therefore can be thought of as a regularization parameter  as well.
However, this does not mean that the product
$\psi(z)\psi(w)$ is divergence-free at arbitrary points $z$ and
$w$. In fact, according to (\ref{fir1}) and (\ref{fir2}),
it is divergent at $z=wq^{\pm 2}$, and therefore we still need
to regularize it by extra means, other than the parameter $q$, to make
its action
well
defined on its modules. This regularization can be achieved through the
normal
ordering. As usual, a product of operators
$\psi_n$ ($\psi_r$) with $n\in\Z$ ($r\in \Z+1/2$) is normal ordered if all
$\psi_n$ ($\psi_r$) with $n<0$ ($r\leq -1/2$)
are moved to the left of all $\psi_n$ ($\psi_r$) with
$n>0$ ($r\geq 1/2$), and in the (R) sector $\psi_0^2$ is normal
ordered in such a way that the normal ordered product
$\psi(z)\psi(z)$ is null, that is, $\psi_0^2=0$. Note that the
same reasoning about the notion of
normal ordering must also be applied to products of the operators
$\al(n)$ and $a$ with $n\in \Z$ and $a\in \C[P]$ to make the action of
products of $X(\ep|z)$ free of divergence. In this case the normal
ordering consists in moving all $\al(n)$ with
$n<0$ to the left of all $\al(n)$ with $n>0$, and $a$ to the
left of $\al(0)$. Therefore, products of $\al(n)$ with $n\in
\Z^*$ in $X(\ep|z)$ as defined by (\ref{XS}) are already normal
ordered, and hence the action of $X(\ep|z)$ on the symmetric
algebra $\Sh$ is well defined by construction.
The remaining relation (\ref{equala2}) is equivalent to
\bac
\{\psi_{n+2},\psi_{m}\}&=&(q^{2}+q^{-2})
\{\psi_{n+1},\psi_{m+1}\},\quad n,m\in \Z,\\
\{\psi_{r+2},\psi_{s}\}&=&(q^{2}+q^{-2})
\{\psi_{r+1},\psi_{s+1}\},\quad r,s\in \Z+{1\over 2},
\ea
and can easily be verified using (\ref{co1}) and (\ref{co2}).
The Clifford algebra with the
anticommutaion relations (\ref{co1}) and
(\ref{co2}) first appeared in the work of Bernard on the
explicit construction of $U_q(so(2n+1))$-standard modules
with level $k=1$ \cite{Ber89}. Then, it was used in \cite{Idzu93,BoWe94}
for the calculation of the $N$-point correlation functions for the
spin-1 XXZ model.

Let  $T^{R}$ ($T^{NS}$), $T^{R}_{even}$
($T^{NS}_{even}$), $T^R_{odd}$ ($T^{NS}_{odd}$)
be the Fock space spanned by the
modes $\{\psi_n,n<0\}$ ($\{\psi_r, r<0\}$),
a subspace of $T^{R}$ ($T^{NS}$) spanned by
an even number of modes $\psi_n$ ($\psi_r$),
and a subspace of $T^R$ ($T^{NS}$) spanned by an
odd number of modes $\psi_n$ ($\psi_r$)
respectively.
Note that in the (R) sector the zero mode $\psi_0$
acts trivially on $T^R$, and  so to make its action
non-trivial we extend the Fock
space $T^R$ by $\C^2$, with basis
\(\{v_+={\scriptsize\left(\begin{array}{l}1\\ 0\end{array}\right)},v_-=
{\scriptsize\left(\begin{array}{l}0\\1\end{array}\right)}\}\) such that
$\psi_n (n\neq 0)$ and $\psi_0$ act as
$\psi_n\otimes 1$ and \( 1\otimes
{\scriptsize\left(\begin{array}{ll}0&1\\1&0\end{array}\right)}\)
on $T^R\otimes \C^2$
respectively.
Putting all the pieces together, we conclude that
$X(\ep|z)$ acts as
\be
X(\ep|z)=S^-_{-\ep}(zq^{-2\ep })S^+_{-\ep }(zq^{2\ep })\otimes
\psi(z)\otimes e^{\ep \al}z^{\ep{\al(0)\over 2}+{1\over 2}}
\ee
on the space $\Sh\otimes (T^R\otimes \C^2)\otimes e^{\al\over 2}
\C[Q]$ in the (R) sector, and on the space
$\Sh\otimes (T^{NS})\otimes \C[Q]$ in the (NS) sector.
Moreover,  the $\uq$-standard modules $V(2\Lambda_0)$,
$V(2\Lambda_1)$ and $V(\Lambda_0+\Lambda_1)$ are isomorphic
to the following subspaces of the latter spaces
\bac
V(2\Lambda_0)&\sim &\Sh\otimes T^{NS}_{even}\otimes \C[2Q]
\oplus\Sh\otimes T^{NS}_{odd}\otimes e^{\al}\C[2Q],\\
V(2\Lambda_1)&\sim &\Sh\otimes T^{NS}_{even}\otimes e^{\al}\C[2Q]
\oplus\Sh\otimes T^{NS}_{odd}\otimes \C[2Q],\\
V(\Lambda_0+\Lambda_1)&\sim &\Sh\otimes
(T^{R}_{even}\otimes v_+\oplus T^{R}_{odd}\otimes v_-)
\otimes e^{\al\over 2}\C[2Q]\\
&&\oplus\Sh\otimes (T^{R}_{odd}\otimes v_+ \oplus
T^{R}_{even}\otimes v_-)\otimes e^{3\al\over 2}\C[2Q],
\ea
and their respective highest weight vectors are given by:
\bac
v_{2\Lambda_0}&=& 1\otimes 1\otimes 1,\\
v_{2\Lambda_1}&=& 1\otimes 1\otimes e^{\al},\\
v_{\Lambda_0+\Lambda_1}&=& 1\otimes 1\otimes v_+\otimes e^{\al\over 2}.
\ea
\\
\noindent Case III:$\quad G(M)$ with arbitrary nonzero $k$:\hfil

In this case, we will not split the ${\cal Z}_q$ algebra into a tensor
product of a group algebra and a new algebra, which is parafermionic
in nature and has been partially described in
\cite{BoLu94}. We will rather construct the $W(M)$ module in terms of
the ${\cal Z}_q$ operators themselves. This will be sufficient
to find an explicit realization of the $\uq$-generalized
Verma module $G(M)=\Sh\otimes W(M)$ since the $\Sh$
is already  constructed in terms of symmetric polynomials of
$\al(n)$.  To this end, we first define
\bea
{\cal Z}(\ep,\epp|z,w)&=&f(\ep,\epp|wz^{-1})
{\cal Z}(\ep|z){\cal Z}(\epp|w),\label{Zepp1}\\
f(\ep,\epp|z)&=&{(q^{-(\ep+\epp)k/2+k-2}z;q^{2k})^
{\ep\epp}_{\infty}\over
(q^{-(\ep+\epp)k/2+k+2}z;q^{2k})^
{\ep\epp}_{\infty}}.
\label{Karim}
\eea
Next, following some ideas in \cite{LeWi84} in the case of
$\LS$ algebra, and which have been well summarized  in \cite{Fodal94}
 in the case of the elliptic algebra,
we introduce the following formal Laurent and power series:
\bea
{\cal Z}(\ep,\epp|z,w)&=&\sum_{n_{1},n_{2}\in\Z}
{\cal Z}(\ep,\epp|n_1,n_2)z^{-n_1}w^{-n_2},\label{expa1}\\
{\cal Z}(\ep|z)&=&\sum_{n\in\Z}
{\cal Z}(\ep|n)z^{-n},\label{expa2}\\
f(\ep,\epp|z)&=&{1\over \sum_{n\geq 0}\~a^{\ep,\epp}_nz^{-n}}
=\sum_{n\geq 0}a^{\ep,\epp}_nz^{-n}.
\label{expa3}
\eea
Relation (\ref{Karim}) and (\ref{expa3}) imply that
\bac
a^{\ep,\epp}_0&=&\~a^{\ep,\epp}_0=1,\\
\sum_{n\geq 0}\~a^{\ep,\epp}_{n}a^{\ep,\epp}_{m-n}&=&
\delta_{m,0},\quad m\geq 0,\\
a^{\ep,\epp}_n&=&\~a^{\ep,\epp}_n=0,\quad n<0.
\ea
Substituting the above expansions back in (\ref{Zepp1}),
and comparing the coefficients of the powers of $wz^{-1}$, we
obtain:
\bea
{\cal Z}(\ep|n_{1}){\cal Z}(\epp|n_2)&=&\sum_{n\geq 0}\~a^{\ep,\epp}_n
{\cal Z}(\ep,\epp|n_1-n,n_2+n),\label{sina1}\\
{\cal Z}(\ep,\epp|n_1,n_2)&=&\sum_{n\geq 0}a^{\ep,\epp}_n
{\cal Z}(\ep|n_1-n){\cal Z}(\epp|n_2+n).
\label{sina2}
\eea
Furthermore, substituting the latter expansions in the quantum generalized
relations (\ref{qZ1}) and (\ref{qZ2}) and using (\ref{sina2}), we arrive at
\bea
{\cal Z}(\ep,-\ep|n_1,n_2)&=&{\cal Z}(-\ep,\ep|n_2,n_1)+Y(\ep|n_1)
\delta_{n_1+n_2,0},
\label{khal1}\\
{\cal Z}(\ep,\ep|n_1,n_2)&=&q^{2\ep}{\cal Z}(\ep,\ep|n_1-1,n_2+1)
+q^{2\ep}{\cal Z}(\ep,\ep|n_2,n_1)-{\cal Z}(\ep,\ep|n_2+1,n_1-1),\nonumber \\
&&\label{khal2}
\eea
respectively, and where
\be
Y(\ep|n)={1\over \q}(q^{kn+\ep\al(0)}-q^{-kn-\ep\al(0)}).
\ee
Relations (\ref{khal1}) and (\ref{khal2}) are useful in the normal
ordering of products of  ${\cal Z}_q$ operators by moving any operator
${\cal Z}(\ep|n_1)$ with $n_1>n_2$ to the right of  ${\cal Z}(\epp|n_2)$,
and the operator ${\cal Z}(+|n)$
to the right of ${\cal Z}(-|n)$. To see this let us examine the
normal ordering of ${\cal Z}(\ep|n_{1}){\cal Z}(\epp|n_2)$ in the
following three nontrivial cases:\\
A.$\quad n_1>n_2$, $\quad \ep=-\epp$:\\
Using relations (\ref{sina1}) and (\ref{khal1}), we obtain
\bac
{\cal Z}(\ep|n_{1}){\cal Z}(-\ep|n_2)&=&
\displaystyle\sum_{n\geq 0}\~a^{\ep,-\ep}_n {\cal Z}(\ep,-
\ep|n_1-n,n_2+n)=\displaystyle\sum_{0\leq n\leq {n_1-n_2-x\over 2}}
\~a^{\ep,-\ep}_n {\cal Z}(-\ep,\ep|n_2+n,n_1-n)\\
&+&\displaystyle\sum_{0\leq {n_1-n_2-x\over 2}< n}
\~a^{\ep,-\ep}_n {\cal Z}(\ep,-\ep|n_1-n,n_2+n)+
\sum_{0\leq n\leq {n_1-n_2-x\over 2}}
\~a^{\ep,-\ep}_nY(\ep|n_1)\delta_{n_1+n_2,0},
\ea
where $x$ is equal to 0 or 1 depending on whether
$n_1-n_2$ is even or odd respectively.
It is therefore clear from the latter relation and
(\ref{sina2}) that the
product ${\cal Z}(\ep|n_{1}){\cal Z}(-\ep|n_2)$ with $n_1>n_2$
can be normal ordered.\\
B. $\quad n_1=n_2$, $\quad\ep=-\epp=+$:\\
Like the previous case, we have
\bac
{\cal Z}(+|n_{1}){\cal Z}(-|n_1)&=&\sum_{n\geq 0}\~a^{+,-}_n {\cal Z}(+,
-|n_1-n,n_1+n),\\
&=&{\cal Z}(-,+|n_1,n_1)+Y(+|0)\delta_{n_1,0}+
\sum_{n>0}
\~a^{+,-}_n {\cal Z}(-\ep,\ep|n_1-n,n_1+n).
\ea
The same argument as in case A holds, and hence the
product ${\cal Z}(+|n_{1}){\cal Z}(-|n_1)$ can be normal ordered as well.\\
C. $\quad n_1>n_2$, $\quad\ep=\epp$:\\
This case is less straightforward. First, relation
(\ref{sina1}) implies that
\be
{\cal Z}(\ep|n_{1}){\cal Z}(\ep|n_2)=\sum_{n\geq 0}\~a^{\ep,\ep}_n
{\cal Z}(\ep,\ep|n_1-n,n_2+n),
\ee
which according to (\ref{sina2}) means that the product
${\cal Z}(\ep|n_{1}){\cal Z}(\ep|n_2)$ can be normal ordered if we can
normal order also any operator ${\cal Z}(\ep,\ep|n_1,n_2)$ with
$n_1>n_2$ by writing
it as a linear combination of operators ${\cal Z}(\ep,\ep|m_1,m_2)$ with
$m_1\leq m_2$.
Relation (\ref{khal2}) allows indeed this second type of normal ordering.
The reason is that
repeated use of this relation leads to
\bac
{\cal Z}(\ep,\ep|n_1+2p,n_1)&=&
q^{2\ep}{\cal Z}(\ep,\ep|n_1,n_1+2p)+q^{2(p-1)\ep}(q^{2\ep}-1)
{\cal Z}(\ep,\ep|n_1+p,n_1+p)\\
&&+\sum_{n=1}^{p-1}
q^{2(n-1)\ep}(q^{4\ep}-1)
{\cal Z}(\ep,\ep|n_1+n,n_1+2p-n),\quad p>0,\\
{\cal Z}(\ep,\ep|n_1+2p+1,n_1)&=&
q^{2\ep}{\cal Z}(\ep,\ep|n_1,n_1+2p+1)\\
&&+\sum_{n=1}^{p}
q^{2(n-1)\ep}(q^{4\ep}-1)
{\cal Z}(\ep,\ep|n_1+n,n_1+2p+1-n),\quad p>0,
\ea
where all the operators
${\cal Z}(\ep,\ep|m_1,m_2)$ in the right hand sides of the above
equations have
$m_1\leq m_2$. Consequently, the product
${\cal Z}(\ep|n_{1}){\cal Z}(\ep|n_2)$ with
$n_1>n_2$ can also be normal ordered.

Since $X(\ep|z)$ acts as $S^-_{-\ep}(zq^{-2\ep })S^+_{-\ep }(zq^{2\ep })
\otimes
{\cal Z}(\ep|z)$  on $\Sh\otimes W(M)\otimes \C[z,z^{-1}]$ it
is clear that $W(M)$ is
spanned by the vectors in the set
\be
\{{\cal Z}(\ep_1|n_1)\dots {\cal Z}(\ep_s|n_s)v_0,\quad
\ep_i=\pm,\quad n_i<0 ,
\quad i=1,\dots,s; \quad s>0\}.
\ee
The condition $n_i<0$ guarantees that the above vectors have
negative degrees as they should (otherwise they are null)
since $\Sh\otimes W(M)$
is a graded $\uq$-highest weight module. Because of the normal
ordering of the ${\cal Z}(\ep|n)$ operators discussed above, the
above spanning set
for $W(G)$ can be reduced further to the smaller set
\bac
H&=&\{{\cal Z}(\ep_1|n_1)\dots {\cal Z}(\ep_s|n_s)v_0,\quad \ep_i=\pm,\quad
n_i\leq n_{i+1};\\
&&\ep_i\leq\ep_{i+1}\quad{\rm if}\quad n_i=n_{i+1} ,
\quad i=1,\dots,s; \quad s>0\},
\ea
where the order $-<+$ is meant. It can easily be seen that the
 set $H$  is a basis for
$W(M)$ since its character coincides with the
one of $W(G)$ as given by (\ref{cha}). Therefore, we have an
explicit construction of the $\uq$-generalized Verma module
$G(M)$ with nonzero level $k$.

\section{Relations in the ${\cal Z}_q$ enveloping algebra}

In this section, we extend the
quantum generalized commutation relations (\ref{qZ1}) and (\ref{qZ2}) to
relations  satisfied
by arbitrary polynomials of $Z(\ep|z)$.
For this purpose
 let us consider the following operators:
\be
{\cal Z}(\ep_{1},\dots,
\ep_{s}|z_{1},\dots,z_{s})=S^-_{\ep_{1}}(z_{1})
\dots S^-_{\ep_{s}}(z_{s})x_{\ep_{1}}(z_{1})
\dots x_{\ep_{s}}(z_{s})S^+_{\ep_{1}}(z_{1})\dots
S^+_{\ep_{s}}(z_{s}),\quad  s>0,
\ee
which are a generalization of (\ref{ZOP}).
They are expressed in terms of the operators
${\cal Z}(\ep|z)$ introduced in (\ref{ZOP}) as:
\bea
\Zep&=&\prod_{1\leq i<j\leq s}
{(q^{-(\ep_{i}+\ep_{j})k/2+k-2}z_{j}z_{i}^{-1};q^{2k})^
{\ep_{i}\ep_{j}}_{\infty}\over
(q^{-(\ep_{i}+\ep_{j})k/2+k+2}z_{j}z_{i}^{-1};q^{2k})^
{\ep_{i}\ep_{j}}_{\infty}}
{\cal Z}(\ep_{1}|z_{1}){\cal Z}(\ep_{2}|z_{2})
\dots {\cal Z}(\ep_{s}|z_{s}),\nonumber\\
&&\label{Zep1}\\
&=&\prod_{2\leq i\leq s}
{(q^{-(\ep_{1}+\ep_{i})k/2+k-2}z_{i}z_{1}^{-1};q^{2k})^
{\ep_{1}\ep_{i}}_{\infty}\over
(q^{-(\ep_{1}+\ep_{i})k/2+k+2}z_{i}z_{1}^{-1};q^{2k})^
{\ep_{1}\ep_{i}}_{\infty}}{\cal Z}(\ep_{1}|z_{1}){\cal Z}(\ep_{2},
\dots \ep_{s}|z_{2},\dots,z_{s}).\nonumber\\
\label{Zep2}
\eea
The above relations can easily be derived from
\bea
S^{+}_{\ep}(z_{1})X_{\epp}(z_{2})&=&
{(q^{-(\ep+\epp)k/2+k+2}z_{2}z_{1}^{-1};q^{2k})^
{\ep\epp}_{\infty}\over
(q^{-(\ep+\epp)k/2+k-2}z_{2}z_{1}^{-1};q^{2k})^
{\ep\epp}_{\infty}}X_{\epp}(z_{2})S^{+}_{\ep}(z_{1}),\\
S^{-}_{\ep}(z_{1})X_{\epp}(z_{2})&=&
{(q^{-(\ep+\epp)k/2+k-2}z_{1}z_{2}^{-1};q^{2k})^
{\ep\epp}_{\infty}\over
(q^{-(\ep+\epp)k/2+k+2}z_{1}z_{2}^{-1};q^{2k})^
{\ep\epp}_{\infty}}X_{\epp}(z_{2})S^{-}_{\ep}(z_{1}).
\label{bo}\eea
Relations (\ref{aZ}) and (\ref{Zep1}) imply that the operators
$\Zep$ commute also with $U_q(\hat{h})$.

Let us now derive the first type of  the quantum
generalized commutation
relations in the ${\cal Z}_q$ enveloping algebra,
which is valid only if $\ep_{r}=-\ep_{r+1}$:
\bac
&& {\cal Z}(\ep_{1},\dots,\ep_{r},\ep_{r+1},\dots ,\ep_{s}
|z_{1},\dots,z_{r},z_{r+1},\dots, z_{s})-
{\cal Z}(\ep_{1},\dots,\ep_{r+1},\ep_{r},\dots, \ep_{s}
|z_{1},\dots,z_{r+1},z_{r},\dots, z_{s})\\
&&=S^{-}_{\ep_1}(z_{1})\dots
S^{-}_{\ep_s}(z_{s})X_{\ep_1}(z_{1})\dots
X_{\ep_{r-1}}(z_{r-1})
[X_{\ep_r}(z_{r}),X_{\ep_{r+1}}(z_{r+1})]\\
&&.X_{\ep_{r+2}}(z_{r+2})\dots
X_{\ep_{s}}(z_{s})S^{+}_{\ep_1}(z_{1})\dots
S^{+}_{\ep_s}(z_{s})\\
&&\displaystyle={\ep_{r}\over
q-q^{-1}} S^{-}_{\ep_1}(z_{1})\dots
S^{-}_{\ep_s}(z_{s})X_{\ep_1}(z_{1})\dots
X_{\ep_{r-1}}(z_{r-1})
\left(\Psi(z_{r+1}q^{\ep_{r} k/2})\delta(z_{r}z_{r+1}^{-1}
q^{-\ep_{r} k})\right.\\
&&\left. -\Phi(z_{r}q^{\ep_{r} k/2})\delta(z_{r}z_{r+1}^{-1}
q^{\ep_{r} k})\right)
X_{\ep_{r+2}}(z_{r+2})\dots
X_{\ep_{s}}(z_{s})
S^{+}_{\ep_1}(z_{1})\dots
S^{+}_{\ep_s}(z_{s})\\
&&=\displaystyle{\ep_{r}\over
q-q^{-1}}{\cal Z}(\ep_{1},\dots\hat{\ep_{r}},\hat{\ep_{r+1}}\dots \ep_{s}
|z_{1},\dots,\hat{z_{r}},\hat{z_{r+1}},\dots z_{s})\\
&&.\displaystyle\left\{\Psi_{0}\delta(z_{r}z_{r+1}^{-1}q^{-\ep_{r} k})
q^{2\sum_{i>r+1}\ep_{i}}
\prod_{i>r+1}\left({1-q^{-2-(\ep_{i}+\ep_{r})k/2}z_{i}z_{
r+1}^{-1}\over 1-q^{2-(\ep_{i}+\ep_{r})k/2}z_{i}z_{
r+1}^{-1}-1}\right)^{\ep_i}\right.\\
&&\displaystyle \left. -\Phi_{0}\delta(z_{r}z_{r+1}^{-1}q^{\ep_{r} k})
q^{-2\sum_{i>r+1}\ep_{i}}\prod_{i<r}
\left({1-q^{-2-(\ep_{i}+\ep_{r})k/2}z_{r+1}z_{
i}^{-1}-1\over q^{2-(\ep_{i}+\ep_{r})k/2}z_{r+1}z_{
i}^{-1}-1}\right)^{\ep_i}\right\},\quad \ep_{r}=\ep_{r+1},
\label{gqZ}\ea
where we have used  (\ref{0PX}),
(\ref{ps}),  (\ref{PSS}), and (\ref{bo}).
Above, we have also introduced the notation
\be
{\cal Z}(\ep_{1},\dots,\hat{\ep_{r}},\hat{\ep_{r+1}}\dots, \ep_{s}
|z_{1},\dots,\hat{z_{r}},\hat{z_{r+1}},\dots, z_{s})=
S^{-}XS^{+},
\ee
with
\bac
S^{\pm}&=&S^{\pm}_{\ep_1}(z_{1})\dots S^{\pm}_{\ep_{r-1}}(z_{r-1})
S^{\pm}_{\ep_{r+2}}(z_{r+2})\dots
S^{\pm}_{\ep_s}(z_{s}),\\
X&=&X_{\ep_1}(z_{1})\dots
X_{\ep_{r-1}}(z_{r-1})X_{\ep_{r+2}}(z_{r+2})\dots
X_{\ep_{s}}(z_{s}),
\ea
where the hat on $\hat{\ep}$ means that the symbol
$\ep$ is omitted.

Using the formal power series
\bac
\displaystyle \left({1-a\over 1-b}\right)^{\ep}&=&
(1-a^{(1+\ep)/2}b^{(1-\ep)/2})
\sum_{n\geq 0}a^{n(1-\ep)/2}b^{n(1+\ep)/2},\\
\prod_{i=1}^{s}(1-z_{i}z^{-1}q^{a_{i}})&=&
\displaystyle \sum_{j_{1},\dots,j_{s}=0,1}(-1)^{\sum_{i=1}^{s}j_i}
q^{\sum_{i=1}^{s}a_i}z_{1}^{j_1}\dots z_{s}^{j_s}z^{-\sum_{i=1}^{s}j_i}.
\ea
we can expand the products in (\ref{gqZ}) as
\bac
&&\displaystyle \prod_{i>r+1}\left({1-q^{-2-(\ep_{i}+\ep_{r})k/2}
z_{i}z_{r+1}^{-1}\over 1-q^{2-(\ep_{i}+\ep_{r})k/2}z_{i}z_{
r+1}^{-1}-1}\right)^{\ep_i}=
\sum_{j_{r+2},\dots,j_{s}=0,1}
\sum_{m_{r+2}\geq j_{r+2},\dots,m_{s}\geq
j_{s}}\\
&&\displaystyle (-1)^{\sum_{i>r+1}j_i}q^{\sum_{i>r+1}
m_i(2\ep_{i}-(\ep_{i}+\ep_{r})k/2)-4\ep_{i}j_{i}}
z_{r+1}^{-\sum_{i>r+1}
m_i}
z_{r+2}^{m_{r+2}}\dots z_{s}^{m_{s}},
\ea
and
\bac
&&\displaystyle \prod_{i<r}
\left({1-q^{-2-(\ep_{i}+\ep_{r})k/2}z_{r+1}z_{
i}^{-1}-1\over q^{2-(\ep_{i}+\ep_{r})k/2}z_{r+1}z_{
i}^{-1}-1}\right)^{\ep_i}=
\sum_{j_{1},\dots,j_{r-1}=0,1}
\sum_{m_{1}\geq j_{1},\dots,m_{r-1}\geq
j_{r-1}}\\
&&\displaystyle (-1)^{\sum_{i<r}j_i}q^{\sum_{i<r}
m_i(2\ep_{i}-(\ep_{i}+\ep_{r})k/2)-4\ep_{i}j_{i}}
z_{1}^{-m_{1}}\dots z_{r-1}^{-m_{r-1}}
z_{r+1}^{\sum_{i<r}
m_i}.
\ea
Substituting these expansions in (\ref{gqZ})
and taking into account the formal Laurent  expansions
\bac
&&{\cal Z}(\ep_{1},\dots, \ep_{s}
|z_{1},\dots, z_{s})
=\sum_{n_{1},\dots, n_{s}\in Z}
{\cal Z}(\ep_{1},\dots, \ep_{s}
|n_{1},\dots, n_{s})z_{1}^{-n_1}\dots
z_{s}^{-n_s},\\
&&{\cal Z}(\ep_{1},\dots\hat{\ep_{r}},\hat{\ep_{r+1}}\dots \ep_{s}
|z_{1},\dots,\hat{z_{r}},\hat{z_{r+1}},\dots z_{s})=\\
&&\displaystyle \sum_{n_{1},\dots,n_{r-1},n_{r+2},\dots, n_{s}\in Z}
{\cal Z}(\ep_{1},\dots,\hat{\ep_{r}},\hat{\ep_{r+1}},\dots, \ep_{s}
|n_{1},\dots,\hat{n_{r}},\hat{n_{r+1}},\dots,
n_{s})z_{1}^{-n_1}\dots z_{r-1}^{-n_{r-1}}
z_{r+2}^{-n_{r+2}}\dots  z_{s}^{-n_s}.
\label{exp}
\ea
we obtain the first quantum generalized
commutation relation satisfied by the Laurent modes:
\bac
&&{\cal Z}(\ep_{1},\dots,\ep_{r},\ep_{r+1},\dots, \ep_{s}
|n_{1},\dots,n_{r},n_{r+1},\dots, n_{s})
-{\cal Z}(\ep_{1},\dots,\ep_{r+1},\ep_{r},\dots, \ep_{s}
|n_{1},\dots,n_{r+1},n_{r},\dots, n_{s})\\
&&\displaystyle =
{\ep_{r}\over
q-q^{-1}}\left(
\sum_{j_{r+2},\dots,j_{s}=0,1}
\sum_{m_{r+2}\geq j_{r+2},\dots,m_{s}\geq
j_{s}}\delta_{n_{r}+n_{r+1},\sum_{i>r+1}m_i}
(-1)^{\sum_{i>r+1}j_i}\right.\\
&&\left. \displaystyle . q^{n_{r}\ep_{r}k+\sum_{i>r+1}
m_i(2\ep_{i}-(\ep_{i}+\ep_{r})k/2)+2\ep_{i}(1-2j_{i})}\right.\\
&&\left. .{\cal Z}(\ep_{1},\dots,\hat{\ep_{r}},\hat{\ep_{r+1}},\dots, \ep_{s}
|n_{1},\dots,n_{r-1},\hat{n_{r}},\hat{n_{r+1}},
n_{r+2}+m_{r+2},\dots, n_{s}+m_{s})\Psi_0\right.\\
&&\left.\displaystyle -\sum_{j_{1},\dots,j_{r-1}=0,1}
\sum_{m_{1}\geq j_{1},\dots,m_{r-1}\geq
j_{r-1}}\delta_{n_{r}+n_{r+1},-\sum_{i<r}m_i}
(-1)^{\sum_{i<r}j_i}\right.\\
&&\left.\displaystyle .q^{-n_{r}\ep_{r}k-2
\sum_{i>r+1}\ep_i+\sum_{i<r}
m_i(2\ep_{i}-(\ep_{i}+\ep_{r})k/2)-4\ep_{i}j_{i}}\right.\\
&&\left. .{\cal Z}(\ep_{1},\dots,
\hat{\ep_{r}},\hat{\ep_{r+1}},\dots, \ep_{s}
|n_{1}+m_{1},\dots,n_{r-1}+m_{r-1},\hat{n_{r}},\hat{n_{r+1}},
n_{r+2},\dots, n_{s})\Phi_0\right),\quad \ep_{r}=-\ep_{r+1}.
\ea
We remark that in the above formula
the following relations are meant:
\bac
\delta_{n_{r}+n_{r+1},\sum_{i>r+1}m_i}&=&
\delta_{n_{r}+n_{r+1},0},\quad {\rm if}\qquad r=1,\\
\delta_{n_{r}+n_{r+1},-\sum_{i<r}m_i}&=&
\delta_{n_{r}+n_{r+1},0},\quad {\rm if}\qquad r=s-1.
\ea
The second type of quantum generalized commutation
relations in the ${\cal Z}_q$ enveloping algebra is derived in the case
$\ep_{r}=
\ep_{r+1}$ as follows:
\bac
&& (z_{r}-z_{r+1}q^{2\ep_{r}}){\cal Z}(\ep_{1},\dots,\ep_{r},
\ep_{r+1},\dots ,\ep_{s}
|z_{1},\dots,z_{r},z_{r+1},\dots, z_{s})\\
&&-(z_{r}q^{2\ep_{r}}-z_{r+1}) {\cal Z}(\ep_{1},\dots,
\ep_{r+1},\ep_{r},\dots, \ep_{s}
|z_{1},\dots,z_{r+1},z_{r},\dots, z_{s})\\
&&=S^{-}_{\ep_1}(z_{1})\dots
S^{-}_{\ep_s}(z_{s})X_{\ep_1}(z_{1})\dots
X_{\ep_{r-1}}(z_{r-1})\\
&&. ((z_{r}-z_{r+1}q^{2\ep_{r}})X_{\ep_r}(z_{r})X_{\ep_{r+1}}(z_{r+1})-
(z_{r}q^{2\ep_{r}}-z_{r+1})X_{\ep_{r+1}}(z_{r+1})X_{\ep_{r}}(z_{r}))\\
&&.X_{\ep_{r+2}}(z_{r+2})\dots
X_{\ep_{s}}(z_{s})
S^{+}_{\ep_1}(z_{1})\dots
S^{+}_{\ep_s}(z_{s})=0,\quad \ep_{r}=\ep_{r+1}.
\ea
Substituting the Laurent expansion (\ref{exp}) in this relation
we obtain this second quantum generalized commutation
relation satisfied by the Laurent modes
\bac
&&{\cal Z}(\ep_{1},\dots,\ep_{r},\ep_{r+1},\dots, \ep_{s}
|n_{1},\dots,n_{r-1},1+n_{r},n_{r+1},\dots, n_{s})\\
&&-q^{2\ep_{r}}{\cal Z}(\ep_{1},\dots,\ep_{r+1},\ep_{r},\dots, \ep_{s}
|n_{1},\dots,n_{r-1},n_{r+1},1+n_{r},\dots, n_{s})\\
&&=q^{2\ep_{r}}{\cal Z}(\ep_{1},\dots,\ep_{r},\ep_{r+1},\dots, \ep_{s}
|n_{1},\dots,n_{r-1},n_{r},1+n_{r+1},\dots, n_{s})\\
&&-{\cal Z}(\ep_{1},\dots,\ep_{r+1},\ep_{r},\dots, \ep_{s}
|n_{1},\dots,n_{r-1},1+n_{r+1},n_{r},\dots, n_{s}),
\quad{\rm if}\quad \ep_{r}=\ep_{r+1}.
\ea

Although the first  quantum generalized
commutation relation looks  completely different from
its classical analogue due to the
appearance of various sums over the indices
$j_i$ and $m_i$ instead of a single sum over a single index
in the classical case \cite{LePr85}, we have checked
that in the limit $q\ra 1$ this quantum generalized
commutation relation  does indeed
reduce to its classical analogue, which in our notation
reads simply as follows:
\bac
&&{\cal Z}(\ep_{1},\dots,\ep_{r},\ep_{r+1},\dots, \ep_{s}
|n_{1},\dots,n_{r},n_{r+1},\dots, n_{s})
-{\cal Z}(\ep_{1},\dots,\ep_{r+1},\ep_{r},\dots, \ep_{s}
|n_{1},\dots,n_{r+1},n_{r},\dots, n_{s})\\
&&=(n_{r}k+2\ep_{r}(\ep_{r+2}+\dots +\ep_{s})+\ep_{r}
h){\cal Z}(\ep_{1},\dots,
\hat{\ep_{r}},\hat{\ep_{r+1}},\dots, \ep_{s}
|n_{1},\dots,n_{r-1},\hat{n_{r}},\hat{n_{r+1}},
n_{r+2},\dots, n_{s}),\\
&& {\rm if}\quad n_{r}+n_{r+1}=0,\quad \ep_r=-\ep_{r+1},\\
&&=2\ep_{r}\sum_{i>r+1}\ep_{i}
{\cal Z}(\ep_{1},\dots,
\hat{\ep_{r}},\hat{\ep_{r+1}},\dots, \ep_{s}
|n_{1},\dots,n_{r-1},\hat{n_{r}},\hat{n_{r+1}},
n_{r+2},\dots,n_{i}+n_{r}+n_{r+1},\dots, n_{s}),\\
&&{\rm if}\quad n_{r}+n_{r+1}>0,\quad \ep_r=-\ep_{r+1},\\
&&=-2\ep_{r}\sum_{i<r}\ep_{i}
{\cal Z}(\ep_{1},\dots,
\hat{\ep_{r}},\hat{\ep_{r+1}},\dots, \ep_{s}
|n_{1},\dots, n_{i}+n_{r}+n_{r+1},\dots,n_{r-1},\hat{n_{r}},\hat{n_{r+1}},
n_{r+2},\dots, n_{s}), \\
&&{\rm if}\quad n_{r}+n_{r+1}<0\quad \ep_r=-\ep_{r+1}.
\ea
The classical analogue of the second quantum generalized
commutation relation is simpler and given by:
\bac
&&{\cal Z}(\ep_{1},\dots,\ep_{r},\ep_{r+1},\dots, \ep_{s}
|n_{1},\dots,n_{r},n_{r+1},\dots, n_{s})\\
&&-{\cal Z}(\ep_{1},\dots,\ep_{r+1},\ep_{r},\dots, \ep_{s}
|n_{1},\dots,n_{r+1},n_{r},\dots, n_{s})=0,
\quad {\rm if}\quad \ep_r=\ep_{r+1}.
\ea
\section{Conclusions}

In this paper, we have introduced a natural quantum
analogue of the ${\cal Z}_q$ algebra with arbitrary level
$k$. In the special cases $k=1$ and $k=2$ our ${\cal Z}_q$ algebra
simplifies considerably and reduces
to the  only presently available results in the literature, which
 become special cases of our general construction. Moreover,
as a new example of a   ${\cal Z}_q$-module, that is,
with $k>2$, we
provide an explicit construction of the basis for
a generalized Verma module. One would like to construct all
${\cal Z}_q$-modules, and from them all $\uq$-modules
 with arbitrary  level $k$, and especially the standard ones.
We believe that the two types of quantum generalized commutation
relations derived in the last section will be useful for the latter
purpose.
Moreover,
one would like to diagonalize the off-critical $ Z_k$ statistical
models   by the
elements of the quantum parafermionic algebra, which is obtained from
the quotient of the ${\cal Z}_q$ algebra by its subalgebra $\C[Q]$, in the
same way that
the off-critical Ising model has been diagonalized
by the quantum Clifford algebra (the special case $k=2$ of the
quantum parafermionic algebra) \cite{Foda94}.

\section*{Acknowledgements}
We are very grateful to O. Foda and G. Watts for
very stimulating discussions and for pointing to us several
references.
A.H.B. is thankful to NSERC for providing him
with a postdoctoral fellowship. The work of L.V.
is supported through funds provided by
NSERC (Canada) and FCAR (Qu\'ebec).

\pagebreak

\end{document}